\documentclass[a4paper,10pt,leqno]{article}
\usepackage{latexsym}
\usepackage[all]{xy}

\usepackage{amssymb} 
\usepackage{amsmath} 
\usepackage{theorem}

\def\Z{{\mathbb{Z}}}
\def\K{{\mathbb{K}}}
\def\A{{\mathcal{A}}}

\DeclareMathOperator{\rank}{rank}

\DeclareMathOperator{\codim}{codim}

\DeclareMathOperator{\Der}{Der}

\numberwithin{equation}{section}

\newcommand{\owari}{\hfill$\square$}

\theoremstyle{break}
\newtheorem{theorem}{Theorem}[section]
\newtheorem{prop}[theorem]{Proposition}
\newtheorem{cor}[theorem]{Corollary}
\newtheorem{lemma}[theorem]{Lemma}
\newtheorem{define}[theorem]{Definition}
\newtheorem{rem}[theorem]{Remark}

\newtheorem{example}[theorem]{Example}

\title{The Euler
multiplicity and addition-deletion theorems 
for multiarrangements}
\author{Takuro Abe\thanks{Supported by 21st Century COE Program 
``Mathematics of Nonlinear Structures via Singularities'' Hokkaido University. } 
\and Hiroaki Terao\thanks{Supported in part by Japan Society for the Promotion
of Science.}
\and Max Wakefield\thanks{Supported by NSF grant \# 0600893 and the NSF Japan program. }}
\date{\today}

\begin{document}

\maketitle

\begin{abstract}
The 
addition-deletion theorems for hyperplane arrangements, 
which were originally shown 
in \cite{T}, provide useful 
ways to construct examples of 
free arrangements. In this article, we prove
addition-deletion theorems for multiarrangements. 
A key to
the generalization is
the definition of a new multiplicity, 
called the {\it Euler
multiplicity},
of a restricted multiarrangement. 
We compute the Euler multiplicities in many cases.
Then we apply
the addition-deletion theorems  
to various arrangements including supersolvable
arrangements and the Coxeter arrangement of type $A_{3} $
to construct free and non-free multiarrangements. 
\end{abstract}
\setcounter{section}{-1}

\section{Introduction}
Let $\A$ be a \textit{hyperplane arrangement}, or simply an
\textit{arrangement}.
In other words,
$\A$  is a 
finite collection of hyperplanes in an $\ell$-dimensional 
vector space $V$ over a field $\K$. 
A \textit{multiarrangement},
which was introduced by Ziegler in \cite{Z},
is a pair $(\A,m)$ consisting of a 
hyperplane arrangement $\A$ and a {\it
multiplicity} $m:\A \rightarrow \Z_{>0}$. 
Define $|m| = \sum_{H\in\A} m(H)$. 
A
multiarrangement $(\A,m)$ 
such that $m(H)=1$ for all $H\in \A$ 
is just a hyperplane arrangement, 
and is sometimes called a 
\textit{simple arrangement}.

Let
$\{x_1,\ldots,x_{\ell}\}$ be a basis for $V^*$. 
Then
$S:=\mbox{Sym}(V^*) \simeq \K[x_1,\ldots,x_{\ell}]$.
When each $H \in \A$ contains the origin, we say that
$\A$ is \textit{central}. Throughout this article,
assume that every arrangement is central.
Let $\Der_{\mathbb K} (S)$ denote 
the set of $\mathbb{K}$-linear derivations
from $S$ to itself. 
For each $H\in \mathcal{A}$ we choose a defining form $\alpha_H$. 
Following Ziegler \cite{Z}, 
we define 
an $S$-module $D(\A,m)$
of a multiarrangement
$(\mathcal{A},m)$ by
$$
D(\mathcal{A},m)=\{ \theta \in \Der_{\mathbb K} (S) \mid \theta (\alpha_H)\in \alpha_H^{m(H)}S
\text{ for all } H\in \mathcal{A}
\}.
$$
If $D(\A, m)$ is a free $S$-module we say that
$(\A, m)$ is a
 {\it free multiarrangement}.  
When $(\A, m)$ is simple, the module coincides with the usual 
module $D(\A)$ of logarithmic derivations
(e.g., \cite[4.1]{OT}).  Thus free multiarrangements 
generalize free arrangements.

When $(\A, m)$ is a free multiarrangement
we define the \textit{exponents} of $(\A, m)$, 
denoted 
by $\exp(\A, m)$, 
to be the multiset of degrees of 
a
homogeneous basis 
$
\{
\theta_1,\ldots,\theta_{\ell}
\}
$ for $D(\A, m)$:
$$
\exp(\A, m):=(\deg (\theta_1), \ldots, \deg(\theta_{\ell})),
$$
where $\deg (\theta_i):=\deg (\theta_i(\alpha))$ for some linear form $\alpha$ 
with 
$\theta_i (\alpha) \neq 0$. 
Then the multiset
$\exp(\A, m)$ does not
depend  upon choice of basis.

In his groundbreaking paper \cite{Z}, Ziegler
writes ``\dots the theory of multiarrangements and their freeness
is not yet in a satisfactory state. In particular,
we do not know any
addition/deletion theorem \dots.''
It is exactly the subject of this article.
Namely, we generalize the \textit{addition-deletion} theorems 
for simple arrangements
\cite{T} 
to multiarrangements in this article. 
Let $(\A, m)$ be a nonempty multiarrangement
and $\ell \geq 2$.
Fix a 
hyperplane $H_0 \in \A$ and let $\alpha_0$ be 
a defining form for $H_{0}.$
To state the addition-deletion theorems for multiarrangements
we need
to define
the deletion
$(\A', m')$
and
the restriction
$(\A'', m'')$.
First, we define the 
{\it deletion} 
as follows:

\begin{define}
\begin{itemize}
\item[(i)] If $m(H_0)=1$, then $\A':= \A \setminus \{H_0\}$ and 
$m'(H)=m(H)$ for all $H \in \A'$.
\item[(ii)] If $m(H_0) \ge 2$, then $\A':= \A$ and for $H \in \A' =\A$, we define
\[
m'(H)=
\left\{
\begin{array}{rl}
m(H) \ \ \ \ \ \  & \mbox{if}\ H \neq H_0 ,\\
m(H_0)-1 & \mbox{if}\ H = H_0.
\end{array}
\right.
\]
\end{itemize}
\end{define}

Next we define the {restriction}
$(\A'', m'')$.   Let
$$\A''=
\{H_0 \cap K \mid K\in \A \setminus \{H_0\}\},$$
which is an arrangement on $H_{0} $. 
We, however, have more than one choice
to define a multiplicity $m''$.
The definition of a suitable
multiplicity $m''$ is crucial.
The canonical definition is probably
\[
m'' (X) = \sum_{{K\in \A \setminus \{H_{0} \}}\atop{K\cap H_{0} = X}} 
m(K),
\]
which is purely combinatorial and
was  used
in \cite{Y1, Y2, Z}
effectively.
In this article, however,
in order to serve our purposes, we introduce a new multiplicity $m^{*}$, 
called the {\it Euler multiplicity}, whose definition is algebraic rather than 
combinatorial.

For
$X\in \A''$ define
\[
\A_{X} = \{H\in \A\mid X\subset H\}
\,\,\,\,
\text{and}
\,\,\,\,
m_{X} = m \mid_{\A_{X} }.
\]
Choose a coordinate system
$(x_{1}, \dots, x_{\ell})$ so that 
$X$ is defined by $x_{1} =x_{2} =0$. 
Let
$\partial_{x_{i} } $ denote
$\displaystyle \frac{\partial}{\partial x_i}$ 
$(1\leq i\leq \ell)$ 
throughout this article.
By Proposition \ref{basiscor} 
we will see that
$D(\A_X,m_X)$ 
has a basis
\begin{equation}
\label{thetapsi}
\theta_X,
\psi_X,
\displaystyle \partial_{x_3}, 
\displaystyle \partial_{x_4}, 
\ldots,
\displaystyle \partial_{x_\ell},
\end{equation} 
such that 
$\theta_{X} \not\in \alpha_{0} \Der_{\mathbb K} (S)$
and
$\psi_{X} \in \alpha_{0} \Der_{\mathbb K} (S)$. 
\begin{define}
The \textit{Euler multiplicity} 
$m^* : \A'' \rightarrow \Z_{>0}$ is defined by 
$
m^*(X):=\deg \theta_X\ (X \in \A'').
$
Then define the {\it restriction} by $(\A'', m^{*})$. 
\label{emultiplicity} 
\end{define}
For $(\A, m)$ and $H_{0} \in\A$ we say the collection
$(\A, m)$, the deletion $(\A', m')$ and
the restriction $(\A'', m^{*})$ is a {\it triple}.   

\noindent
\begin{rem} 
When $(\A, m)$ is simple
the Euler derivation can be chosen as $\theta_X$. 
In this case,
$m^* \equiv 1$, so 
$(\A'', m^*)$ is simple.
\end{rem}

For 
$\theta \in D(\A, m)$
define 
$\overline{\theta}\in D(\A'')$
by 
$
\overline{\theta}(\overline{f}):=\overline{\theta(f)}
$
for
 $\overline{f} \in \overline{S}:=S/\alpha_0 S$, 
where 
$\overline{f}$ is the image of an element $f \in S$ by the canonical 
projection $S \rightarrow \overline{S}$. 
In Proposition \ref{exact}
we obtain an exact sequence
\[
0
\longrightarrow
D(\A', m')
\stackrel{\alpha_{0}\cdot }{\longrightarrow}
D(\A, m)
\stackrel{\pi}{\longrightarrow}
D(\A'', m^{*} ),
\]
 where $\alpha_{0}\cdot $ denotes the multiplication by
$\alpha_{0}$ and $\pi(\theta) = \overline{\theta}.$

Roughly speaking, 
the addition-deletion theorems
state that the freeness of any two of the triple,
under a condition concerning their exponents,
imply the freeness of the third. 
The following four addition-deletion theorems are 
the multiarrangement versions of
Theorems 4.46 (1),
4.49,
4.46 (2),
and
4.50
 in \cite{OT}.
The ideas behind the proofs are very similar to 
those in \cite{OT}.  However, because of the
indispensability of the  Euler multiplicity,
we include the proofs.

\begin{theorem}
If $(\A,m)$ and
$(\A',m')$ are both free,
then 
there exists a basis $\{\theta_1,\ldots,\theta_{\ell}\}$ for 
$D(\A',m')$ such that,
for some $k \in \{1,\ldots, \ell\}$,
$\{\theta_1,\ldots,\theta_{k-1}, \alpha_0 \theta_k,
\theta_{k+1}, \ldots,\theta_{\ell}\}$ is 
a basis for $D(\A,m)$.
\label{basis}
\end{theorem}

\begin{theorem}[Deletion]
Assume that $(\A,m)$ and $(\A'',m^*)$ are both free and 
$\exp(\A'',m^*) \subset \exp(\A,m)$. Then $(\A',m')$ is also 
free.
\label{deletion}
\end{theorem}

\begin{theorem}[Restriction]
Assume that $(\A,m)$ and $(\A',m')$ are both free. Take a basis 
$\{\theta_1,\ldots,\theta_k,\ldots, \theta_{\ell}\}$ for $D(\A',m')$ 
as in Theorem \ref{basis}. 
Then $\{\overline{\theta_1},\ldots,
\overline{\theta_{k-1}},
\overline{\theta_{k+1}},
\ldots,
\overline{\theta_{\ell}}\}
$ is a basis for $D(\A'',m^*)$.
\label{basis2}
\end{theorem}

\begin{theorem}[Addition]
Assume that $(\A',m')$ and $(\A'',m^*)$ are both free and 
$\exp(\A'',m^*) \subset \exp(\A',m')$. Then 
$(\A,m)$ is also free.
\label{addition}
\end{theorem}

Summarizing these results we follow Cartier \cite{Cartier}
to 
obtain the following addition-deletion theorem for 
multiarrangements.

\begin{theorem}[Addition-Deletion]
Let $(\A,m)$ be a nonempty multiarrangement in an $\ell$-dimensional 
vector space $V$, $H_0 \in \A$ 
and 
let $(\A,m),(\A',m'),\ (\A'',m^{*})$ be the triple with respect to $H_0$. Then any two of the following 
statements imply the third:
\begin{itemize}
\item[(i)]
$(\A,m)$ is free with $\exp(\A,m)=(d_1,\ldots,d_{\ell}).$
\item[(ii)]
$(\A',m')$ is free with $\exp(\A',m')=(d_1,\ldots,d_{\ell}-1).$
\item[(iii)]
$(\A'',m^{*})$ is free with $\exp(\A'',m^*)=(d_1,\ldots,d_{\ell-1}).$
\end{itemize}
\label{addition-deletion}
\end{theorem}

Applying Addition
Theorem \ref{addition} repeatedly, we can
inductively
construct the following class of
free multiarrangements.  

\begin{define}
The class $\mathcal{IFM}$ of 
\textit{inductively free  multiarrangements} is the 
smallest class of multiarrangements 
which satisfies the following two conditions.
\begin{itemize}
\item[(1)] 
The empty arrangement $\emptyset_{\ell}$ in an
$\ell$-dimensional vector space is contained in $\mathcal{IFM}$ for 
$\ell \ge 0$.
\item[(2)] 
For a multiarrangement $(\A,m)$, if there exists $H \in \A$ such that 
$(\A',m') \in \mathcal{IFM},\ 
(\A'',m^*) \in \mathcal{IFM}$, and 
$\exp(\A',m') \supset \exp(\A'',m^*)$, then 
$(\A,m) \in \mathcal{IFM}$. 
\end{itemize}
\label{IFM}
\end{define}

\begin{rem}
The intersection of the class of $\mathcal{IFM}$
with the class of simple arrangements
is equal to the class of 
inductively free arrangements,
$\mathcal{IF}$
\cite[Definition 4.53]{OT}.
\end{rem}

The outline of this article is as follows. 
In Section 1, we introduce some definitions and 
recall some known results 
in arrangement theory which will be used later. 
In Section 2, we prove Theorem \ref{basis} 
and Deletion Theorem \ref{deletion}. 
In Section 3, we prove 
Restriction
Theorem \ref{basis2} and Addition Theorem \ref{addition}. 
In Section 4, we compute explicit values of 
the Euler multiplicities in many cases.
Applying the addition-deletion theorems together with the
computations, in Section 5, we find
multiplicities $m$ such that the multiarrangement $(\A, m)$ 
is free for various arrangements $\A$ 
including supersolvable arrangements
and the
Coxeter arrangement
of type $A_{3}$.

\section{Preliminaries}
In this section we fix some notation and 
introduce some results about multiarrangements which will be used later. For 
hyperplane arrangement theory, 
we refer the reader to \cite{OT}. 
For a multiarrangement $(\A,m)$, define  
$$
Q(\A,m):=\prod_{H \in \A} \alpha_H^{m(H)}.
$$

The $S$-module 
$\Der_{\K}(S)$
of  $\K$-linear $S$-derivations
has the natural basis:
$$
\Der_{\K}(S)
= \bigoplus_{i=1}^{\ell} S\displaystyle 
\partial_{x_i}.
$$ 
We say a nonzero element $\theta
=\sum_{i=1}^{\ell} f_i \partial_{x_i}
 \in \Der_{\K}(S)$ is \textit{homogeneous of degree $p$} if 
$f_i$ is zero or a homogeneous polynomial of degree $p$ in $S$ for 
$1 \le i \le \ell$.  
Recall the $S$-submodule 
$$
D(\mathcal{A},m) = \{\theta \in \Der_{\K}(S)\ |\ \theta(\alpha_H) \in S\cdot \alpha_H^{m(H)}\ (\forall 
H \in \mathcal{A})\}
$$
of $\Der_{\K} (S)$ and 
a multiarrangement $(\mathcal{A},m)$ is 
\textit{free} if 
$D(\mathcal{A},m)$ is a free $S$-module. 
The fact that
the module
$D(\mathcal{A},m)$ is reflexive
(e.g., see Theorem 5 in \cite{Z})
implies
the following proposition.

\begin{prop}
\label{ranktwo}
A multiarrangement $(\mathcal{A},m)$ is free for any
multiplicity $m$ whenever
$r(\A)
:=
\codim_{V}(\bigcap_{H\in\A} H) \leq 2$. 
\end{prop}

For $\theta_1,\ldots,\theta_{\ell} \in D(\A,m)$, we define the 
$(\ell \times \ell)$-matrix 
$M(\theta_1,\ldots,\theta_{\ell})$ as the matrix whose $(i,j)$-entry 
is $\theta_j(x_i)$. In general, it is difficult to determine whether a given 
multiarrangement is free or not. 
However, 
using the following criterion
(see Theorem 8 in \cite{Z} and Theorem 4.19 in \cite{OT}), 
we can verify that a candidate for a basis
is actually a basis.

\begin{theorem}[Saito-Ziegler's criterion]
Let
$\theta_1,\ldots,\theta_{\ell}$ be derivations in $D(\A,m)$. Then 
$\{\theta_1,\ldots,\theta_{\ell}\}$ forms a basis for $D(\A,m)$ if and only if 
$$
\det M(\theta_1,\ldots,\theta_{\ell}) \in \K^{*} \cdot Q(\A,m).
$$ 
In particular, if $\theta_1,\ldots,\theta_{\ell}$ are all homogeneous, 
then $\{\theta_1,\ldots,\theta_{\ell}\}$ 
forms a basis for $D(\A,m)$ if and only if 
the following two conditions are satisfied:
\begin{itemize}
\item[(i)]
$\theta_1,\ldots,\theta_{\ell}$ are independent over $S$.
\item[(ii)] 
$\sum_{i=1}^{\ell} \deg (\theta_i)=\sum_{H \in \A} m(H).$
\end{itemize}
\label{Saito}
\end{theorem}

Let $V_i$ be vector spaces over $\K$ and 
$(\mathcal{A}_i, m_i)$ be multiarrangements in $V_i
\,\,\, (i=1,2)$. Let us define 
their product
$
(\mathcal{A}_1
\times
\mathcal{A}_2
, 
m_{1} \times m_{2})
$
in the vector space 
$
V_1 \oplus V_2$ by the following manner:
\begin{gather*}
\mathcal{A}_1\times 
\mathcal{A}_2
:=
\{
H_1 \oplus V_2
\mid
 H_1 \in \mathcal{A}_1\}
\cup
\{
V_1 \oplus H_2
\mid
 H_2 \in \mathcal{A}_2\},\\
(m_{1} \times m_{2})(H_1 \oplus V_2):=m_1(H_1),\\
(m_{1} \times m_{2})(V_1 \oplus H_2):=m_2(H_2).
\end{gather*}
The following lemma is a special case of Lemma 1.4
in \cite{ATW}.

\begin{lemma}
$$
D(\mathcal{A}_1 \times \A_2, m_1 \times m_2)\simeq S \cdot D(\mathcal{A}_1,m_1) \oplus 
S\cdot D(\mathcal{A}_2,m_2),
$$
where $S = {\rm Sym}
((V_1 \oplus V_2)^{*})$. 
\label{key1}
\end{lemma}

We will use the following 
lemma 
in this article repeatedly.
For the proof see
\cite[Theorem 4.42]{OT} for example.

\begin{lemma}
Let $M = \oplus_{i=0}^{\infty} M_{i}  $
be a free graded $S$-module with a homogeneous basis
$\eta_1,\ldots, \eta_\ell$.
Suppose $\deg\eta_{i} = {d_{i} }
\,\,
(1 \le i \le \ell)$  with
$d_{1} \leq \dots \leq d_{\ell}$.
Assume that there exist elements
$\theta_1,\ldots, \theta_k\ 
(1 \le k \le \ell)
$ in $M$ 
which satisfy the following two 
conditions:
\begin{itemize}
\item[(i)]
$\deg (\theta_i)=d_i\ (i=1,\ldots,k)$.
\item[(ii)]
$\theta_i \not \in S \theta_1 + S \theta_2 + \ldots + S \theta_{i-1}\ (i=1,\ldots,k).$
\end{itemize}
Then $\theta_1,\ldots,\theta_k$ can be extended to a basis for $M$.
\label{basis extend}
\end{lemma}

\section{Deletion}
In this section we prove Theorem \ref{basis} and Deletion Theorem \ref{deletion}. 
We use the notation 
$(d_1,\ldots,d_{\ell})_{\le}$ to indicate 
$d_1 \le \cdots \le d_\ell$.

\medskip

\noindent
\textbf{Proof of Theorem \ref{basis}.} 
Let $(d_1,\ldots,d_{\ell})_{\le}$ be the exponents of $(\A',m')$ 
and 
$(d_1,\ldots,d_{k-1},e_k,\ldots,e_{\ell})_{\le}$ 
be the exponents of 
$(\A,m)$ 
such that $e_k \neq d_k$. 
Choose a basis $\{\theta_1,\ldots,\theta_{k-1},\psi_k,
\ldots,\psi_{\ell}\}$ for $D(\A,m)$ with $\deg(\theta_i)=d_i$ and 
$\deg (\psi_i)=e_i$. 
Because $\theta_1,\ldots,\theta_{k-1}$ are contained in 
$D(\A',m')$ and satisfy the two 
conditions 
in Lemma \ref{basis extend}, we can find a basis 
$\{\theta_1,\ldots,\theta_{k-1},\theta_k,\ldots,\theta_{\ell}\}$ for $D(\A',m')$ 
with $\deg (\theta_i)=d_i\ (k \le i \le \ell)$. 
Since $\alpha_0 \theta_k \in D(\A,m)$, 
$$
\alpha_0 \theta_k =\sum_{i=1}^{k-1} a_i \theta_i +
\sum_{i=k}^{\ell} b_i \psi_i\ (a_i,b_i \in S).
$$
Given that $\theta_1,\ldots,\theta_k$ are independent over $S$, there exists some $j,\ 
j \ge k$ such that $b_j \neq 0$. Hence, 
$$
\deg (\alpha_0 \theta_k)=d_k+1 \ge \deg(\psi_j) \ge \deg(\psi_k) = e_k.
$$
Moreover, since $\psi_k \in D(\A',m')$, 
$$
\psi_k =\sum_{i=1}^{k-1} a_i \theta_i +
\sum_{i=k}^{\ell} b_i \theta_i\ (a_i,b_i \in S).
$$
A similar argument as the above implies 
$$
\deg (\psi_k)=e_k \ge \deg(\theta_j) \ge \deg(\theta_k) = d_k.
$$
The assumption that $e_k \neq d_k$ implies that $e_k=d_k+1$. 
Noting that $\deg(\alpha_0 \theta_k)=d_k+1=e_k$, Lemma \ref{basis extend} shows that 
the elements 
$\{\theta_1,\ldots, \theta_{k-1},\alpha_0 \theta_k\}$, which are contained in 
$D(\A,m)$, can be extended to a basis $\{\theta_1,\ldots, \theta_{k-1},
\alpha_0 \theta_k,\theta_{k+1}',\ldots,\theta_{\ell}'\}$ for 
$D(\A,m)$. Then Theorem \ref{Saito} implies 
$\{\theta_1,\ldots, \theta_{k-1},
\theta_k,\theta_{k+1}',\ldots,\theta_{\ell}'\}$ is a basis for 
$D(\A',m')$. \owari
\medskip

Let $(\A, m)$ be a multiarrangement and $H_{0} \in \A$.
Recall the restriction
$$\A''=
\{H_0 \cap K \mid K\in \A \setminus \{H_0\}\},$$
which is an arrangement on $H_{0} $. 
Let $X\in \A''$.  
Note that $(\A_X,m_X)$ can be decomposed into a direct product of 
a multiarrangement in $\K^2$ and the empty arrangement in 
$X\simeq \K^{\ell-2}$.
Choose a coordinate system 
$(x_{1}, \dots , x_{\ell})$ so that 
$\alpha_{0} = x_{1} $ and 
$X = \{x_{1} =x_{2} = 0\}$.

\begin{prop}
\label{basiscor} 
We may choose a basis
\[
\theta_{X}, \psi_{X}, 
\partial_{x_{3} },
\dots,
\partial_{x_{\ell} }
\]
for $D(\A_{X}, m_{X})$ such that 
$\theta_{X} 
\not\in\alpha_{0} \Der_{\K}(S) $ 
and
$\psi_{X} \in\alpha_{0}  \Der_{\K}(S)$. 
\end{prop}

\noindent
{\bf Proof.}
Let
$
(\A'_{X}, m'_{X}) 
$
be the deletion of 
$ 
(\A_{X}, m_{X})
$
with respect to
$
H_{0}
$.
Then
$
(\A_{X}, 
m_{X}) $ 
and 
$
(\A'_{X},
m'_{X} 
)
$ 
are both
free 
by Proposition \ref{ranktwo}.
It follows from Lemma \ref{key1} and Theorem \ref{basis}
that there exists a
homogeneous
basis 
$\{\theta_{1}, 
\theta_{2}, 
\partial_{x_{3}},
,
\dots
,
\partial_{x_{\ell}}
\}$ 
for 
$D(\A'_{X}, m'_{X} )$
such that
$
\{
\theta_{1} , 
x_{1}  \theta_{2} , 
\partial_{x_{3} },
,
\dots
,
\partial_{x_{\ell} }
\}$ 
 is a basis for 
$D(\A_{X}, 
m_{X})$. 
Define
$
\theta_{X} :=
\theta_{1}$
and
$
\psi_{X}
:= 
x_{1} 
 \theta_{2}.
$ 
It suffices to show that
$\theta_{1} \not\in  
x_{1}  \Der_{\K}(S)$.
If
$
\theta_{X}
=
\theta_{1} 
\in  
x_{1}  \Der_{\K}(S)$
then
$
\theta'_{1} 
:=
\theta_{1} /x_{1} 
\in  
D(\A'_{X}, m'_{X} )$.
This contradicts
the assumption that
$\{\theta_{1} 
=
x_{1} 
\theta'_{1}, \theta_{2} , 
\partial_{x_{3} },
,
\dots
,
\partial_{x_{\ell} }
\}$ is a basis
for 
$D(\A'_{X}, m'_{X} )$.
\owari

\medskip
\noindent
Using the
derivation
$\theta_{X}$ in Proposition \ref{basiscor} we may define
the Euler multiplicity
\[
m^{*} (X)
=
\deg \theta_{X} 
\]
 as in Definition \ref{emultiplicity}.

In Section 0 we defined
the map
$\pi:
D(\A,m) \rightarrow D(\A'')$ 
by $\pi(\theta) = \overline{\theta}$
for $\theta \in D(\A, m)$. 
Note that $\pi$ is well-defined because
$
\overline{\theta(f)} 
=
\overline{\theta(g)} 
$ 
if $f-g\in\alpha_{0} S$. 
Let $\alpha_{0} \cdot : D(\A', m') \rightarrow D(\A, m)$ 
be the multiplication map by $\alpha_{0} $.

\begin{prop} 
We have an
exact sequence
$$
0 \longrightarrow D(\A',m')
\stackrel{\alpha_0\cdot}{\longrightarrow} D(\A,m) 
\stackrel{\pi}{\longrightarrow} 
D(\A'',m^*).$$
\label{exact} 
\end{prop} 

\noindent
\textbf{Proof.} The injectivity of $\alpha_{0} \cdot$ 
and 
the exactness at $D(\A, m)$ are both obvious.
So it suffices to show that
$\pi(\theta)$ lies in
$D(\A'',m^*)$ for $\theta\in
D(\A,m)$. Let $X\in\A''$.  
Note that $
D(\A, m)\subseteq
D(\A_{X}, 
m_{X})$.
We use the notation in the proof
of Proposition \ref{basiscor}.
Moreover, by Lemma \ref{basis extend}, we may assume
$\theta_{X} (x_{i} ) \in \K[x_{1} , x_{2} ]$ $(i = 1, 2)$. 
Thus we obtain
$
\theta_{X} (x_{2} ) \in
(x_{1} , x_{2}^{m^{*}(X)}  ) S
$,
or equivalently
$
\pi(\theta_{X}) (\overline{x_{2}} )
\in
\overline{x_{2}}^{m^{*}(X )} \overline{S}.
$ 
Because
$\pi(\psi_{X} ) = 0$ and
$\pi(\partial_{x_{i} } ) (\overline{x_{2}} )=0
$
for $3\leq i\leq \ell$,
we have $\pi(\theta)(\overline{x_{2}} )
\in
\overline{x_{2}}^{m^{*}(X )} \overline{S}$ 
for all $\theta\in D(\A, m)$. 
\owari

\medskip

To show Deletion Theorem \ref{deletion}, we need the following lemma.

\begin{lemma}
Let $(\A'',m^*)$ be a free multiarrangement with exponents 
$(d_1,\ldots,d_{\ell-1})_{\le} $. Assume the elements 
$\theta_1,\ldots,\theta_k\ (1 \le k \le \ell-1)$ 
in $D(\A,m)$ satisfy the following two conditions:
\begin{itemize}
\item[(i)]
$\deg (\theta_i)=d_i\ (i=1,\ldots,k-1).$
\item[(ii)]
$\deg (\theta_k) < d_k$.
\end{itemize}
Then there exists $p,\ 1 \le p \le k$ such that 
\begin{equation}
\theta_p \in S\theta_1+\ldots + S \theta_{p-1} + \alpha_0 D(\A',m').
\label{star}
\end{equation}
\label{basis3}
\end{lemma}

\noindent 
\textbf{Proof.} 
Assume that for all $i,\ 1 \le i \le k$, condition (\ref{star}) is not true. 
Then $\overline{\theta}_1,\ldots,\overline{\theta}_{k-1}$ satisfy the two conditions 
in Lemma \ref{basis extend}. So 
$\overline{\theta}_1,\ldots,\overline{\theta}_{k-1}$ 
can be extended to a basis for $
D(\A'',m^*)$. 
Since $\deg (\overline{\theta_k}) < d_k$, 
$$
\overline{\theta_k}=\sum_{i=1}^{k-1} \overline{a_i} \overline{\theta_i}
$$
for $a_i \in {S}$. This implies
$
\theta_k \in S\theta_1+\ldots + S \theta_{k-1} + \alpha_0 D(\A',m'),
$
which is a contradiction. \owari
\medskip

\noindent 
\textbf{Proof of Deletion Theorem \ref{deletion}.} 
Put
\begin{eqnarray*}
\exp(\A,m)&=&(d_1,\ldots,d_{\ell})_\le,\\
\exp(\A'',m^*)&=&(d_1,\ldots,d_{k-1},d_{k+1},\ldots,d_{\ell})_\le.
\end{eqnarray*}
We may assume that $d_k< d_{k+1}$ or $k= \ell$. First assume 
$d_k < d_{k+1}$. 
Take a basis $\{\theta_1,\ldots,\theta_{\ell}\}$ for $D(\A,m)$ with $\deg(\theta_i)=d_i$. 
Since $\deg (\theta_k)=d_k < d_{k+1}$, Lemma \ref{basis3} shows that 
there exists some $p,\ 1 \le p \le k$ such that 
$$
\theta_p \in S \theta_1 +\cdots +S \theta_{p-1} + \alpha_0 D(\A',m').
$$
Hence we may assume that $\theta_p \in \alpha_0 D(\A',m')$. Then Theorem \ref{Saito} implies 
$\{\theta_1,\ldots,\theta_{p-1},\theta_p/\alpha_0,\theta_{p+1},\ldots,\theta_{\ell}\}$ is a basis 
for $D(\A',m')$. Next assume $k= \ell$, then 
$\exp(\A'',m^*)=(d_1,\ldots,d_{\ell-1})$. 
If $\overline{\theta_i} \in \bar{S} \bar{\theta}_1+
\cdots +\bar{S} \bar{\theta}_{i-1}$ for some $
i,\ 1 \le i \le \ell-1$, then we can use the same argument as above. 
If $\overline{\theta_i} \not \in \bar{S} \bar{\theta}_1+
\cdots +\bar{S} \bar{\theta}_{i-1}$ for all $
i,\ 1 \le i \le \ell-1$, then Lemma \ref{basis extend} shows 
$\{\overline{\theta_1},\ldots,\overline{\theta_{\ell-1}}\}$ is a basis 
for $D(\A'',m^*)$. Hence 
$$
\theta_{\ell} \in 
S \theta_1 +\cdots +S \theta_{\ell-1} + \alpha_0 D(\A',m'),
$$
and the same argument as above completes the proof. \owari
\medskip

\section{Addition and Restriction}
In this section we prove Restriction Theorem \ref{basis2} and 
Addition Theorem \ref{addition}. 
First, for each $X \in \A''$, let us fix a hyperplane $H_X \in \A \setminus 
\{H_0\}$ 
such that $\overline{H_X}:=H_0 \cap H_X=X$. 
Let $m_0$ denote $m(H_0)$. Recall the definition of 
$\theta_X,\psi_X \in D(\A_X,m_X)$ in Proposition
$\ref{basiscor}$. 
Denote 
\begin{eqnarray*}
e_X:=\deg(\theta_X)\ \mbox{and}\ d_X:=\deg(\psi_X).
\end{eqnarray*}

\begin{lemma}
Let $(\A,m)$ be a multiarrangement in $\K^2$ with exponents $(d,e)$. 
Fix a line $H_0=\{\alpha_0=0\} \in \A$. By Theorem \ref{basis}, there exists a 
basis $\{\theta,\psi\}$ for $D(\A,m)$ such that  
$\deg (\theta)=e,\ \deg(\psi)=d$ and that $\overline{\theta} \neq 0,\ 
\overline{\psi} =0$. Then  
$d-m_0 \ge 0$.
\label{additional}
\end{lemma}

\noindent
\textbf{Proof.} 
We may assume that $S \simeq \K[x_1,x_2]$ and 
$\alpha_0=x_1$. 
If 
$\psi(x_1) = 0$, 
then 
$\theta(x_1) \neq 0$ and Theorem \ref{Saito} implies 
$Q(\A,m)\in\K^{*} \cdot\theta(x_1)\psi(x_2)$. Also we have 
$x_1|\psi(x_2)$ and $x_1^{m_0}|\theta(x_1)$. This implies
$x_1^{m_0+1}|Q(\A,m)$, which is a contradiction. 
So we may assume that 
$\psi(x_1)\neq
0$. 
Therefore,
$x_1^{m_0}|\psi(x_1)$
and thus
$\deg (\psi)=d \ge m_0$. 
\owari
\medskip

\begin{prop}
For all $X \in \A''$, we have $d_X - m_0 \ge 0$.
\label{polynomial}
\end{prop}

\noindent 
\textbf{Proof.} 
Since $(\A_X,m_X)$ can be decomposed into a direct product of 
a multiarrangement in $\K^2$ and the empty arrangement in 
$X\simeq\K^{\ell-2}$,  Lemma \ref{key1} and 
Lemma \ref{additional} complete the proof. \owari
\medskip

By Proposition \ref{polynomial}, 
we make the following key definition.

\begin{define}
Define a polynomial $B=B(\A'',m^*)$ by
$$
B(\A'',m^*):=\alpha_0^{m_0-1} \prod_{X \in \A''} \alpha_{H_X}^{d_X - m_0}.
$$
\end{define}

\begin{lemma}
For any $\theta \in D(\A',m')$, we have 
$\theta(\alpha_0) \in (\alpha_0^{m_0},B(\A'',m^*))$.
\label{ideal}
\end{lemma}

\noindent 
\textbf{Proof.} 
Take any $X \in \A''$ and consider the $S$-module 
$D(\A_X',m_X')$, which contains $D(\A',m')$ as a submodule. 
Since $X$ is of codimension two and $\ell \ge 2$, 
$(\A_X',m_X')$ is free with exponents 
$(e_X,d_X-1,0,\ldots,0)$. 
By Proposition \ref{basiscor}, we have basis elements 
$\theta_X$ and $\psi_X'$ for $D(\A'_X,m'_X)$ such that 
$\deg (\theta_X)=e_X,\ \deg (\psi_X')=d_X-1$ and that $\theta_X$ and $\alpha_0 \psi_X'$ are 
basis elements for $D(\A_X,m_X)$. First we show $D(\A_X',m_X')\alpha_0
:=\{\theta(\alpha_0) |\ 
\theta \in D(\A_X',m_X')\} \subseteq (\alpha_0^{m_0}, 
\alpha_0^{m_0-1} \alpha_{H_X}^{d_X-m_0})$. 
We may assume that $\alpha_0=x_1$ and 
$\alpha_{H_X}=x_2$. Then
$
\{\theta_X,\psi_X',\partial_{x_3},
\ldots,\partial_{x_{\ell}}\}
$
is a basis for $D(\A_X',m_X')$. Since $\theta_X (x_1) \in x_1^{m_0}S$ 
and $\partial_{x_i}(x_1)=0\ (2 \le i \le \ell)$, it suffices to show 
$\psi_X'(x_1) \in (x_1^{m_0}, 
x_1^{m_0-1} x_2^{d_X-m_0})$. We may assume that 
$\psi_X'$ is a derivation of $\K[x_1,x_2]$ by Lemma \ref{key1}. 
Thus there exist $f \in \K[x_1,x_2]$ and $g \in \K[x_2]$ such that
$$
\psi_X'(x_1)=x_1^{m_0-1}(x_1 f(x_1,x_2)+g(x_2)).
$$
Note that $\deg(\psi_X'(x_1))=d_X-1$, so $\deg(g(x_2))=
d_X-m_0$. Hence 
$$
\psi_X'(x_1) \in (x_1^{m_0},x_1^{m_0-1} x_2^{d_X-m_0}).
$$
So we have
\begin{eqnarray*}
D(\A',m')\alpha_0 \subseteq \bigcap_{X \in \A''} D(\A_X',m_X')\alpha_0
&\subseteq& \bigcap_{X \in \A''} (\alpha_0^{m_0},\alpha_0^{m_0-1}\alpha_{H_X}^{d_X-m_0})\\
&=&(\alpha_0^{m_0}, \alpha_0^{m_0-1} \prod_{X \in \A''} \alpha_{H_X}^{d_X-m_0}). 
\ \ \ \ \ \ \ \ \ \ \ \ \ \ \ \ \ \ \ \ \ \ \ \ \ \ \ \ \ \ \ \ \ \ \ \ \ \ \ 
\square
\end{eqnarray*}
\medskip

\noindent 
\textbf{Proof of Restriction Theorem \ref{basis2}.} 
Recall that we have a basis 
$\{\theta_1,\ldots,\theta_k,\ldots,\theta_\ell\}$ for 
$D(\A',m')$ such that 
$\{\theta_1,\ldots,\alpha_0\theta_k,\ldots,\theta_\ell\}$ is a basis for 
$D(\A,m)$. Noting that 
$$
|m|
=\sum_{H \in \A}m(H)
=m_0+ 
\sum_{X \in \A''} (e_X+d_X-m_0),
$$
we have 
$$
\deg(B(\A'',m^*))=m_0-1+\sum_{X \in \A''} (d_X -m_0)
=|m|-1-\sum_{X \in \A''}e_X=|m'|-|m^*|.
$$
Assume that $\deg(\theta_k)<\deg(B(\A'',m^*))$. Then 
Lemma \ref{ideal} implies that $\theta_k(\alpha_0) \in \alpha_0^{m_0}S$. This is equivalent to 
$\theta_k \in D(\A,m)$, which contradicts Theorem \ref{basis}. 
Hence, $\deg(\theta_k) \ge \deg (B(\A'',m^*)).$ 
This inequality implies  
\begin{eqnarray*}
\sum_{i \neq k} \deg (\theta_i)&=&
|m'|-\deg (\theta_k)\\
&\le& |m'| -\deg(B(\A'',m^*)) 
= |m'|-(|m'|-|m^*|)
=|m^*|.
\end{eqnarray*}
To complete the proof by using Theorem \ref{Saito}, it suffices to show 
$\{\overline{\theta_1},\ldots,\overline{\theta_{k-1}},
\overline{\theta_{k+1}},\ldots, \overline{\theta_{\ell}}\}$ is independent over 
$\overline{S}$. Assume that there exist $a_i \in S\ (i=1,\ldots,k-1,k+1,\ldots, \ell)$ such that 
$\sum_{i \neq k} \overline{a_i}\overline{\theta_i}=0$. This implies that 
there exists some $\theta \in \Der_{\K}(S)$ 
such that 
$$
\sum_{i \neq k} a_i \theta_i = \alpha_0 \theta.
$$
Since $\theta_1,\ldots,\theta_{k-1},\theta_{k+1},\ldots,\theta_{\ell}$ lie in $D(\A,m)$, 
we can see that $\theta \in D(\A',m')$, and this implies $\overline{a_i}=0$ for all $i$. \owari
\medskip

\noindent 
\textbf{Proof of Addition Theorem \ref{addition}.} 
Denote 
\begin{eqnarray*}
\exp(\A',m')&=&(d_1,\ldots,d_{\ell})_\le,\\
\exp(\A'',m^*)&=&(d_1,\ldots,d_{k-1},d_{k+1},\ldots, d_{\ell})_\le . 
\end{eqnarray*}
Choose a basis $\{\theta_1,\ldots,\theta_{\ell}\}$ for $D(\A',m')$ such that 
$\deg (\theta_i)=d_i\ (i=1,\ldots, \ell)$. 
We may assume that $d_k < d_{k+1}$ or $k= \ell$. Note that 
$\deg (B(\A'',m^*))=|m'|-|m^*|=d_k$ in this case. Hence, 
Lemma \ref{ideal} implies any $\theta_j$ satisfying $\deg (\theta_j) < \deg (B(\A'',m^*))=d_k$
is contained in $D(\A,m)$. 
First assume that $d_k < d_{k+1}$. 
Consider the following condition:
\begin{equation}
\mbox{For all}\ \theta_j\ \mbox{with}\ \deg(\theta_j)=d_k,\ 
\mbox{it holds that}\ \theta_j \in D(\A,m). \label{assume}
\end{equation}
If (\ref{assume}) is true, then 
$\theta_1,\ldots,\theta_k \in D(\A,m)$. Applying 
Lemma \ref{basis3}, we can see that there exists $p,\ 1 \le p \le k$ such that 
$$
\theta_p \in S \theta_1 +\ldots+S \theta_{p-1} + \alpha_0 D(\A',m').
$$
Thus we may assume that 
$\theta_p \in \alpha_0 D(\A',m')$.
This implies that $\alpha_0^{m_0}|\alpha_0 \cdot
\det M(\theta_1,\ldots,\theta_p/\alpha_0,\ldots,\theta_\ell)$, 
which is a contradiction. 
So, there exists some $p,\ 1 \le p \le k$ such that 
$\deg (\theta_p)=d_k$ and $\theta_p \not \in D(\A,m)$. 
Let us put
$$
\theta_p(\alpha_0)=f_p \alpha_0^{m_0} +c_p B(\A'',m^*)
$$
for some $f_p,c_p \in S$. 
Since $\deg(\theta_p)=d_k=\deg (B(\A'',m^*))$ and $\theta_p \not \in D(\A,m)$, we may assume that 
$c_p=1$. Similarly, for $j \neq p$, put 
$$
\theta_j(\alpha_0)=f_j \alpha_0^{m_0} +c_j B(\A'',m^*)
$$
for some $f_j,c_j \in S$. Define $\eta_j:=\theta_j - c_j \theta_p\ (j \neq p)$ and 
$\eta_p:=\alpha_0 \theta_p$. Then 
$\eta_1,\ldots,\eta_{\ell} \in D(\A,m)$. Theorem \ref{Saito} implies 
$\{\eta_1,\ldots,\eta_{\ell}\}$ is a basis for $D(\A,m)$. 

Next assume that $k= \ell$. If (\ref{assume}) is true, then 
$\theta_1,\ldots,\theta_{\ell} \in D(\A,m)$, which is a contradiction. 
Hence there exists some $p,\ 1 \le p \le \ell$ such that 
$\deg (\theta_p)=d_{\ell}=\deg (B(\A'',m^*))$ and $\theta_p \not \in D(\A,m)$. Then the same argument as above 
completes the proof. \owari

\section{Euler multiplicities}\label{emulti}

To apply the addition-deletion theorems the  
computation of the Euler multiplicities $m^*$ of the restriction is  
crucial. In general, computing the Euler multiplicities is difficult. On the other hand, using results from \cite{W} and \cite{WY} we can compute the Euler multiplicities in the following cases.

\begin{prop}\label{combemult}

Let $X\in \mathcal{A}''$ where  $\mathcal{A}''$ is the 
restriction to $H_0\in \mathcal{A}$ and $m_0=m(H_0)$. Suppose $k=|\mathcal{A}_X|$ and $m_1=\mathrm{max}\{ m(H)| H\in \mathcal{A}_X\backslash \{H_0\} \}$.

\begin{itemize}
\item[(1)] If $k=2$ then $m^*(X)=m_1$.

\item[(2)] If $2m_0\geq |m_X|$ then 
$m^*(X)=|m_X|-m_0$.

\item[(3)] If $2m_1\geq |m_X|-1$ then $m^*(X)=m_1$.

\item[(4)] If $|m_X|\leq 2k-1$ and $m_0>1$ then $m^*(X)=k-1.$

\item[(5)] If $|m_X| \leq 2k-2$ and $m_0=1$ then $m^*(X)=|m_X|-k+1.$

\item[(6)] If $m_X\equiv 2$ then $m^*(X)=k$.

\item[(7)] If $k=3$, $2m_0\leq |m_X|$, and $2m_1\leq |m_X|$ then $m^*(X) 
=\left\lfloor \frac{|m_X|}{2}\right\rfloor$.

\end{itemize}
\end{prop}

\noindent 
{\bf Proof.} Without loss of generality we can assume $\ell =2$, $H_0=\{ x_1=0\}$, $Q(\mathcal{A}_X)=x_1x_2Q$, and $Q(\mathcal{A}_X,m_X)=x_1^{m_0}x_2^{m_1}\tilde{Q}$ for some $Q,\tilde{Q}\in S$.  Then  
for case (1) we have $\theta_X=x_2^{m_1}\partial_{x_2}$ and $\psi_X=x_1^{m_0} 
\partial_{x_1}$ in the notation of Proposition \ref{basiscor}. Thus, for case (1) we have $m^*(X)=m_1$. In case (2), if $2m_0 
\geq |m_X|$ then the fact that $m^*(X)=|m_X|-m_0$ follows from \cite{WY} because the smallest degree  
derivation is of the form $\theta=\frac{Q(\mathcal{A}_X,m_X)}{x_1^{m_0}}\partial_{x_2}$. 
Case (3) is similar to case (2). The only difference is that 
one of a basis element 
is of the form $\theta=
\frac{Q(\mathcal{A}_X,m_X)}{x_2^{m_1}}\partial_{x_1}$ and 
$\overline{\theta}=0$. 
Now, suppose that 
$ |m_X|\leq 2k-1$. Then the exponents are
$(|m_X|-k+1, k-1)$ and
$\varphi =\frac{Q(\mathcal{A}_X,m_X)}{Q(\mathcal{A}_X)}
(x_1\partial_{x_1}+x_2\partial_{x_2})$ 
can be chosen as a basis element
by \cite{WY}.  
In case (4), $\varphi$ is divisible by $x_1$ 
and hence  
$m^*(X)=k-1$. 
In case (5), $\varphi$ is not divisible  
by $x_1$.
Since $\varphi$ 
is a basis element of the smallest degree, we have
$m^*(X)=\deg\varphi=|m_X|-k+1$. 
In case (6), if  
$m_X\equiv 2$ then the exponents are $(k,k)$ 
(see Proposition 5.4 in \cite{SoT}).  
In case (7) the  
formula given by Wakamiko
in \cite{W} for the smallest degree generator is not  
divisible by $x_1$. Thus, in case (7) $m^*(X)=\left\lfloor 
\frac{|m_X|}{2}
\right\rfloor$. 
\owari

\

The next example shows that even when the exponents are combinatorially determined the Euler multiplicities may depend on the position of the hyperplanes. 

\begin{example}\label{hmm}
Consider the class of two-dimensional multiarrangements 
$({\mathcal A}_\xi ,m)$ given by the defining polynomial $\tilde{Q}_{\xi}=x_1^4x_2^3(x_1-x_2)(x_1-\xi x_2)$ 
where $\xi \in \mathbb{K}-\{ 0,1\}$. Then a basis for  
$D(\mathcal{A}_\xi ,m)$ for all $\xi \in \mathbb{K}-\{ 0,1\}$ is the following derivations $$\theta_1=x_1^4 
\partial_{x_1} +\left[ (1+\xi (1+\xi))x_1x_2^3-\xi (1+\xi )x_2^4\right] 
\partial_{x_2}$$ and $$\theta_2=x_2^3(x_1-x_2)(x_1-\xi x_2)\partial_{x_2}.$$ Suppose that $({\mathcal A}_\xi ,m)$ is of the form $(\mathcal{A}_X,m_X)$ for some $(\mathcal{A},m)$ and $X\in \mathcal{A}''$ where $H_0=\{ x_1=0\}$. Then the basis $\{ \theta_1,\theta_2\}$ shows that $$m^*(X)=\left\{ \begin{array}{ll}
5& \text{ if } \xi =-1\\
4& \text{ otherwise}
\end{array}_.\right.$$

\end{example}


\section{Applications}

In this section, we apply 
the addition-deletion theorems 
together with the computations of the Euler multiplicities in Proposition \ref{combemult} to construct free and non-free multiarrangements. 

\begin{define}
Let $\A=\{H_1,\ldots,H_n\}$ be a simple arrangement. 
Then $[ m_1,\ldots ,m_n]\in \mathbb{Z}^n_{>0}$ is a {\it free multiplicity} for $\mathcal{A}$ if $(\mathcal{A},m)$ is a free multiarrangement where $m(H_i)=m_i$ for all $1\leq i\leq n$.
\label{freemultiplicity}
\end{define}

It is difficult 
to determine which multiplicities are free for a fixed simple arrangement. 
At least 
the 
following 
proposition provides an infinite number of free multiplicities for an arbitrary 
free arrangement.

\begin{prop}
Let $\A$ be a free simple arrangement with 
$\exp(\A)=(1,d_2,\ldots,d_{\ell})$. 
Fix one hyperplane $H_0 \in \A$ and consider a 
multiarrangement $(\A,m)$ where $m$ is defined by 
\[
m(H)=
\left\{
\begin{array}{rl}
1\  & \mbox{if}\ H \neq H_0 ,\\
m_0 & \mbox{if}\ H = H_0.
\end{array}
\right.
\]
Then $(\A,m)$ is free with 
$
\exp(\A,m)=(m_0,d_2,\ldots,d_{\ell})
$.
\label{onemultiplicity}
\end{prop}

\noindent
\textbf{Proof}. 
Let $(\A,m),(\A',m')$ and 
$(\A'',m^*)$ be the triple with respect to $H_0$. 
Recall the restricted multiarrangement $(\A'',m'')$, where 
$m''(X)=|\A_X|-1$ for all $X \in \A''$ which is 
defined by Ziegler in \cite{Z}.  
It is proved in \cite{Z} that if $\A$ is free with 
$\exp(\A)=(1,d_2,\ldots,d_{\ell})$, then 
$(\A'',m'')$ is also free with 
$\exp(\A'',m'')=(d_2,\ldots,d_{\ell})$. 
Let $X \in \A''$. By Proposition \ref{combemult} (2) and (4), 
$m^*(X)=|\A_X|-1=m''(X)$. 
To finish the proof, apply Addition Theorem \ref{addition}. \owari
\medskip

In the next example, we exhibit a free multiarrangement that is 
not inductively free by using Proposition \ref{onemultiplicity}. 

\begin{example}
Recall the arrangement $\A$ in Example 4.59, based on a pentagon, 
in \cite{OT}, which is due to 
K. Brandt and J. Keaty. This arrangement is free with  
exponents $(1,5,5)$, but it is 
not inductively free. Fix $H_0 \in \A$ which is not 
the infinite hyperplane. Then 
by Proposition \ref{onemultiplicity}, the multiplicity $m$ defined by 
\[
m(H)=
\left\{
\begin{array}{rl}
1 & \mbox{if}\ H \neq H_0 ,\\
2 & \mbox{if}\ H = H_0.
\end{array}
\right.
\]
is a free multiplicity of $\A$ and 
$\exp(\A,m)=(2,5,5)$. Since $\A$ is not inductively free, 
to show $(\A,m)$ is not 
inductively free, it suffices to show 
that any deletion $(\A',m')$ with respect to $H \in \A \setminus \{H_0\}$ 
is not free. 
By Proposition \ref{combemult} (1), (3) and (5), 
the restricted multiarrangement 
$(\A'',m^*)$ with respect to $H$ has Euler multiplicity 
$m^*=[2,1,1,1,1]$. Hence  
$\exp(\A'',m^*)=(2,4)$. Now Deletion Theorem \ref{deletion} implies that 
$(\A',m')$ is not free, so $(\A,m)$ is not inductively free. 
\label{nonmultifree}
\end{example}



\begin{define}
An arrangement $\A$ is \textit{totally free} (or \textit{totally 
non-free}) if $(\A,m)$ is free (respectively non-free) for 
any multiplicity $m$ on $\A$.
\label{total}
\end{define}

\begin{rem}
If $ \ell \le 2$, then any arrangement is totally free
by Proposition \ref{ranktwo}. 
Also, if $\A_1$ and $\A_2$ are both totally free, then so is 
$\A_1 \times \A_2$ by 
Lemma \ref{key1}. Consequently, 
any Boolean arrangement is totally free.
\label{coning}
\end{rem}

\begin{example}
Let $\A$ be an arrangement consisting of four generic hyperplanes in 
$\K^{3}$. 
Let 
\[
Q(\A, m) = x_1^{a} x_2^{b} x_3^{c} (x_1+x_2+x_3)^{d} 
\]
with $1\leq a\leq b\leq c\leq d$. 
We will show that $\A$ is a totally non-free arrangement.
Suppose that $(\A, m)$ is free with minimum $|m|$.
Let 
$\exp(\A, m) = (d_{1}, d_{2}, d_{3})_{\leq} $.
Let $H_{0} = \{x_1=0\}$ and
$\exp(\A'', m^{*}) = (e_{1}, e_{2})_{\leq}$.
If $a=1$, then $(\A', m')$ is Boolean with exponents
  $(b, c, d)$.  Thus $(e_{1}, e_{2}) \subset (b, c, d)$.
This is a contradiction
because $e_{1} + e_{2} = b+c+d$.  So we may assume $2\leq a$.

Case 1. If $d_{1} < e_{1} $, then $(\A', m')$ is free, which is a 
contradiction.

Case 2. If $d_{1} = e_{1} $ and $d_{2} \leq e_{2} $, 
then $(\A', m')$ is free, which is a contradiction.

Case 3. If $d_{1} = e_{1} $ and $d_{2}  > e_{2} $, 
then 
this is a contradiction because 
$a+b+c+d=d_{1} +d_{2} +d_{3} 
> e_{1} + e_{2} + d_{3} 
=
b+c+d+d_{3} \geq a+b+c+d$. 

Case 4. If $d_{1} > e_{1} $ and $b+c\leq d$, 
then $d_{1} > e_{1} = b+c \geq a+b$.
This is a contradiction because 
$x_1^{a} x_2^{b} (\partial_{x_1} - \partial_{x_2})\in D(\A, m)$.

Case 5. If $d_{1} > e_{1} $ and $b+c > d$, 
then 
$d_{1} > e_{1} $ 
and
$d_{1} \geq e_{1}+1\geq e_{2}  $.
This is a contradiction because 
$a+b+c+d=d_{1} +d_{2} +d_{3} 
> e_{1} + e_{2} + d_{3} 
=
b+c+d+d_{3} \geq a+b+c+d$.

\label{TNF}
\end{example}

\begin{rem}
In general, an arrangement can be neither 
totally free nor totally non-free (see Example 14 in \cite{Z}). 
Also note that 
the example by Edelman and Reiner in \cite{ER2} 
is a non-free simple arrangement which admits 
a free multiplicity. 
\label{EdelmanReiner}
\end{rem}

Let us consider supersolvable arrangements defined by 
Stanley in \cite{S}. (The following definition is equivalent to the original 
definition.)

\begin{define}
An arrangement $\A$ is \textit{supersolvable} if 
there exists a filtration
$$
\A=\A_{r} \supset \A_{r-1} \supset \cdots \supset \A_2 \supset \A_1
$$
such that 
\begin{itemize}
\item[(1)] 
$\rank(\A_i)=i\ (i=1,\ldots,r).$
\item[(2)] 
For any $H,H' \in \A_i$, there exists some $H'' \in \A_{i-1}$ such that 
$H \cap H' \subset H''$.
\end{itemize}
\label{supersolvable}
\end{define}

\begin{rem}
It is shown in \cite{terao86} that an arrangement is supersolvable if and only if it is 
\textit{fiber type}. 
\end{rem}

Let us consider a multiarrangement $(\A,m)$ for a supersolvable arrangement $\A$. 
It is shown in \cite{JT} and \cite{St2} that $m \equiv 1$ is a free multiplicity. 
The following theorem gives another sufficient condition 
for $m$ to be a free multiplicity.

\begin{theorem}
Let $(\A,m)$ be a multiarrangement such that $\A$ is supersolvable with a filtration 
$\A=\A_{r} \supset \A_{r-1} \supset \cdots \supset \A_2 \supset \A_1$ and 
$r \ge 2$. Let 
$m_i$ denote the multiplicity $m|_{\A_i}$ and 
$\exp(\A_2,m_2)=(d_1,d_2,0,\ldots,0)$. 
Assume that for each $H' \in \A_d \setminus \A_{d-1}$, 
$H'' \in \A_{d-1}\ (d=3,\ldots,r)$ and $X:=H' \cap H''$, 
it holds that 
\begin{equation}
\A_X=\{H',H''\} \label{eq:000}
\end{equation}
or that
\begin{equation}
m(H'') \ge \sum_{X \subset H \in (\A_d \setminus \A_{d-1})} m(H) -1. \label{eq:001}
\end{equation}
Then $(\A,m)$ is inductively
free with $$\exp(\A,m)
=(d_1,d_2,|m_3|-|m_2|,\ldots,|m_r|-|m_{r-1}|,0,\ldots,0).$$
\label{freeparallel}
\end{theorem}

\noindent
\textbf{Proof.} 
Let us put $d_i:=|m_i|-|m_{i-1}|\ (i=3,\ldots,r)$. 
We may assume that 
$$
\left\{\prod_{i=1}^d x_i=0 \right\} \subseteq \A_d
$$
for all $d,\ 1 \le d \le \ell$. 
We prove by an induction on $r$. When $r=2$, there is nothing to prove. 
Assume $r \ge 3$ and $(\A_{r-1},m_{r-1})$ is free with  
$\exp(\A_{r-1},m_{r-1})=(d_1,d_2,d_3,\ldots,d_{r-1},0,\ldots,0)$. We show 
that $(\A_r,m_r)$ is free with $\exp(\A_r,m_r)=
(d_1,d_2,d_3,\ldots,d_{r-1},d_{r},0,\ldots,0)$. 
Let $H_r \in \A_r \setminus \A_{r-1}$ 
and $(\A_r'',m_r^*)$ be the restricted multiarrangement with 
respect to $H_r$. Since $\A$ is supersolvable, 
$\A_r''=\A_{r-1}|_{H_r}$. Also, the conditions 
(\ref{eq:000}), (\ref{eq:001}),
Proposition \ref{combemult} (1) and (3) imply 
$m_r^*(X)=m_{r-1}(H)$ where
$H \in \A_{r-1}$ and $X=H \cap H_r$. Hence 
$(\A_r'',m_r^*)$ is free with 
$\exp(\A_r'',m_r^*)=(d_1,d_2,\ldots,d_{r-1},0,\ldots,0).$ 
To complete the proof, apply this argument and 
Addition Theorem \ref{addition} repeatedly. \owari

\medskip

Theorem \ref{freeparallel} gives many free multiplicities on 
supersolvable arrangements. For the remainder of this article, 
assume that $\ell =3$ and 
we consider the following supersolvable multiarrangement $(\A,m)$.

\begin{define}
The \textit{Coxeter multiarrangement of type $A_3$} 
can be defined by the following polynomial:
$$
Q(\A,m)=x_1^{m_1} (x_1-x_3)^{m_2} (x_1-x_2)^{m_3} x_2^{m_4} (x_2-x_3)^{m_5} x_3^{m_6}.
$$
\label{a3}
\end{define}

The filtration is given by 
\begin{eqnarray*}
\A_1:&=&\{x_1=0\},\\
\A_2:&=&\{x_1x_2(x_1-x_2)=0\},\\
\A_3:&=&\{x_1x_2x_3(x_1-x_2)(x_1-x_3)(x_2-x_3)=0\}.
\end{eqnarray*}
It is shown in \cite{Sa1}, \cite{Sa2} and 
\cite{T4} that $m=[m_1\,\ldots,m_6]=[m,m,m,m,m,m]$ for $m \in \Z_{>0}$ 
is a free multiplicity 
of $\A$. 
Now 
we obtain the following corollary by using Theorem \ref{freeparallel}.

\begin{cor}
Let $(\A,m)$ be the Coxeter multiarrangement of type $A_3$. 
Assume that $m_1 \ge \max\{m_3,m_4\},\ m_1 \ge m_2+m_6-1,\ 
m_4 \ge m_5+m_6-1$ and $m_3 \ge m_2+m_5-1$. Then 
$(\A,m)$ is inductively
free with
\begin{multline*}
\exp(\A,m)\\
=
\left\{
\begin{array}{ll}
(\left\lfloor \frac{m_1+m_3+m_4}{2}\right\rfloor,
\lceil \frac{m_1+m_3+m_4}{2} \rceil,m_2+m_5+m_6)
 & \mbox{if~~}\ m_1 \le m_3+m_4-1 ,\\
(m_1,m_3+m_4,m_2+m_5+m_6) & \mbox{if~~}\ m_1 > m_3+m_4-1.
\end{array}
\right.
\end{multline*} 
\end{cor}

Remark that for any $H \in \A$
the Euler multiplicity on $H$ can be calculated by 
Proposition \ref{combemult} (1), (2), (3) and (7).

\begin{example}
Let $\A$ be a Coxeter arrangement of type $A_3$. Then 
the multiplicity $m:=
[1,1,2,2,1,1]$ is free by Theorem \ref{freeparallel}. 
However, the multiplicity $k:=[2,1,1,1,2,1]$ is not free. 
Assume $k$ is free. It is shown in \cite{Sa1} and 
\cite{Sa2} that $\A$ is free with $\exp(\A)=(1,2,3)$. Also, Proposition 
\ref{onemultiplicity} implies $m_0:=[2,1,1,1,1,1]$ is free with 
$\exp(\A,m_0)=(2,2,3)$. 
Then Theorem \ref{basis} implies $\exp(\A,k)=(2,2,4)$ or $(2,3,3)$. 
However, for the restricted multiarrangement $(\A'',k^*)$  
with respect to $x_2-x_3=0$, we can see that 
$k^*=[2,2,2]$. Hence $\exp(\A'',k^*)=(3,3)$, which contradicts 
Restriction Theorem \ref{basis2}. Hence 
$k=[2,1,1,1,2,1]$ is not a free multiplicity of $\A$. 
We note that the non-freeness criterion in \cite{ATW} also 
shows that $k$ is not a free multiplicity.
\label{a3nonfree}
\end{example}

 \noindent
 Takuro Abe\\
 Department of Mathematics, \\
 Hokkaido University, \\
 Sapporo 060-0810, Japan. \\
 abetaku@math.sci.hokudai.ac.jp

 \vspace{3mm}

 \noindent
 Hiroaki Terao\\
 Department of Mathematics, \\
 Hokkaido University, \\
 Sapporo 060-0810, Japan. \\
 terao@math.sci.hokudai.ac.jp

 \vspace{3mm}

 \noindent
 Max Wakefield\\
 Department of Mathematics, \\
 Hokkaido University, \\
 Sapporo 060-0810, Japan. \\
 wakefield@math.sci.hokudai.ac.jp


\begin{thebibliography}{ATW}

 \bibitem[ATW]{ATW} T. Abe, H. Terao and M. Wakefield, The characteristic polynomial of a
 multiarrangement. 
 \textit{Advances in Math.} \textbf{215} (2007), 825--838.

\bibitem[C]{Cartier} P. Cartier,
Les arrangements d'hyperplans:
un chapitre de g\'eometrie combinatroire.
In \textit{S\'eminaire Bourbaki} 1980/81.
Lecture Notes in Math.
{\bf 901}, Springer Verlag, 1981, pp.1--22.


 \bibitem[ER]{ER2} P. H. Edelman and V. Reiner, 
 A counterexample to Orlik's conjecture.
 \textit{Proc. Amer. Math. Soc.} \textbf{118} (1993), 927--929.


\bibitem[JT]{JT} {M. Jambu and H. Terao},
Free arrangements of hyperplanes and supersolvable 
lattices. 
\textit{Advances in Math.} {\bf 52}
(1984),
248--258.

 \bibitem[OT]{OT} P. Orlik and H. Terao, \textit{Arrangements of hyperplanes.}
 Grundlehren der Mathematischen Wissenschaften, 
\textbf{300}. Springer-Verlag, Berlin, 1992.

\bibitem[Sa1]{Sa1} {K. Saito},
On the uniformization of complements of 
discriminant loci. In: 
\textit{Conference Notes. Amer. Math. Soc. Summer 
Institute, Williamstown} (1975). 

\bibitem[Sa2]{Sa2} {K. Saito},
Theory of logarithmic differential forms and logarithmic vector fields. 
\textit{J. Fac. Sci. Univ. Tokyo Sect. IA  Math}. 
\textbf{27} (1980), 265--291. 


\bibitem[SoT]{SoT} L. Solomon and H. Terao, 
The double Coxeter arrangement.
\textit{Comment. Math. Helv.}
64 (1998), 237--258.

 \bibitem[St1]{S} R. P. Stanley, Supersolvable lattices. \textit{Algebra Universalis} 2 (1972),
 197--217. 

 \bibitem[St2]{St2} R. P. Stanley, $T$-free arrangements of hyperplanes.
 In: \textit{Progress in Graph Theory}. Academic Press (1984), p. 539.



 \bibitem[T1]{T} H. Terao, Arrangements of hyperplanes and their freeness I, II. \textit{J. Fac. Sci. Univ. Tokyo Sect. IA  Math}. \textbf{27} 
(1980), 293--320. 



\bibitem[T2]{terao81} {H. Terao},
Generalized exponents of a free arrangement of hyperplanes and
Shephard-Todd-Brieskorn formula. \textit{Invent. math.} 
\textbf{63}  (1981),
159--179.

\bibitem[T3]{terao86} {H. Terao},
Modular elements of lattices and topological fibration.
\textit{Advances in Math.} {\bf 62}
(1986),
135--154.


 \bibitem[T4]{T4} {H. Terao}, Multiderivations of Coxeter arrangements. 
\textit{Invent. math.} 
\textbf{148} (2002), 659--674. 

\bibitem[Waka]{W} A. Wakamiko, On the Exponents of $2$-Multiarrangements. 
\textit{Tokyo. J. Math. }
\textbf{30} (2007), 99--116.



\bibitem[WY]{WY} M. Wakefield and S. Yuzvinsky, Derivations of an effective divisor 
on the complex projective line. 
\textit{Trans. Amer. Math. Soc.}
\textbf{359} (2007), 4389--4403.

 \bibitem[Y1]{Y1} M. Yoshinaga, Characterization of a free arrangement and
 conjecture of
 Edelman and Reiner. \textit{Invent. math.} \textbf{157} (2004), no.2,
 449--454.

 \bibitem[Y2]{Y2} M. Yoshinaga, On the freeness of 3-arrangements. 
\textit{Bull. London Math. Soc.} \textbf{37} (2005), no. 1, 126--134. 

 \bibitem[Z]{Z}
 G. M. Ziegler,
 Multiarrangements of hyperplanes and their freeness. in
 {\it Singularities} (Iowa City,
 IA, 1986), 345--359, Contemp. Math., {\bf 90}, Amer. Math. Soc.,
 Providence, RI, 1989.

 \end{thebibliography}
\end{document}